\documentclass[12pt]{article}
\author{Michaela Cully-Hugill\\ School of Science\\ The University of New South Wales Canberra, Australia\\m.cully-hugill@unsw.edu.au \and  Tim Trudgian\footnote{Supported by Australian Research Council Future Fellowship FT160100094.} \\
School of Science\\ The University of New South Wales Canberra, Australia \\
  t.trudgian@adfa.edu.au}
  \title{Two explicit divisor sums}
\usepackage{indentfirst}
\usepackage{url}
\usepackage{enumerate}
\usepackage{amsthm}
\usepackage{amsmath}
\usepackage{mathtools}
\usepackage{comment}
\usepackage{fullpage}
\usepackage{amssymb}
\usepackage{booktabs}
\usepackage{siunitx}
\usepackage{xcolor}
\usepackage{graphicx}

\newtheorem{thm}{Theorem}

\newtheorem{cor}{Corollary}

\begin{document}
\maketitle
\begin{abstract}
\noindent 
We give explicit bounds on sums of $d(n)^{2}$ and $d_{4}(n)$ where $d(n)$ is the number of divisors of $n$ and $d_{4}(n)$ is the number of ways of writing $n$ as a product of four numbers. In doing so we make  a slight improvement on the upper bound for class numbers of quartic number fields. \end{abstract}
\section{Introduction}\label{intro}
\noindent
Let $d(n)$ denote the number of divisors of $n$, and for $k\geq 2$ let $d_{k}(n)$ denote the number of ways of writing $n$ as a product of $k$ integers. Using the following Dirichlet series
\begin{equation}\label{s1}
\sum_{n=1}^{\infty} \frac{d(n)^{2}}{n^{s}} = \frac{\zeta^{4}(s)}{\zeta(2s)}, \quad \sum_{n=1}^{\infty} \frac{d_{k}(n)}{n^{s}} = \zeta^{k}(s), \quad (\Re(s) >1),
\end{equation}
we have, via Perron's formula, that
\begin{equation}\label{s2}
\sum_{n\leq x} d_{k}(n) \sim \frac{1}{(k-1)!} x (\log x)^{k-1},
\end{equation}
and
\begin{equation}\label{s2a}
\sum_{n\leq x} d(n)^{2} \sim \frac{1}{\pi^{2}} x (\log x)^{3}.
\end{equation}
The purpose of this article is to consider good explicit versions of (\ref{s2a}) and (\ref{s2}) when $k=4$. The rolled-gold example is $\sum_{n\leq x} d(n)$, which is gives a bound of the form (\ref{s2}) for $k=2$. Berkan\'{e}, Bordell\`{e}s and Ramar\'{e} \cite[Thm. 1.1]{Berk} gave several pairs of values $(\alpha, x_{0})$ such that 
\begin{equation}\label{study}
\sum_{n\leq x} d(n) = x(\log x + 2 \gamma -1) + \Delta(x),
\end{equation}
holds with $|\Delta(x)| \leq \alpha x^{1/2}$ for $x\geq x_0$, and where $\gamma$ is Euler's constant. One such pair given, which we shall use frequently, is $\alpha=0.397$ and $x_0= 5560$. 

The best known bound for (\ref{study}) is $\Delta(x) = O(x^{131/416 + \epsilon})$ by Huxley \cite{Huxley}. It seems hopeless to give a bound on the implied constant in this estimate. Therefore weaker\footnote{It is worth mentioning Theorem 1.2 in (\cite{Berk}), which gives $|\Delta(x)|\leq 0.764 x^{1/3} \log x$ for $x\geq 9995$. This could be used in our present approach, but the improvement is only apparent for very large $x$.}, yet-still-explicit bounds such as those in \cite{Berk} are very useful in applications.


Bordell\`{e}s \cite{Border} considered the sum in (\ref{s2}), and showed that when $k\geq 2$ and $x\geq1$ we have
\begin{equation}\label{hall}
\sum_{n\leq x} d_{k}(n) \leq x \left( \log x + \gamma + \frac{1}{x}\right)^{k-1}.
\end{equation}
This misses the asymptotic bound in (\ref{s2}) by a factor of $1/(k-1)!$. Nicolas and Tenenbaum (see \cite{Border}, pg. 2) were able to meet the asymptotic bound in (\ref{s2}), by showing that for a fixed $k$, and $x\geq 1$, we have
\begin{equation}\label{lounge}
\sum_{n\leq x} d_{k}(n) \leq \frac{x}{(k-1)!} \left( \log x +(k-1)\right)^{k-1}.
\end{equation}
We improve these results in our first theorem.
\begin{thm}\label{parlour}
For $x\geq 2$ we have 
\begin{equation}\label{chaise}
\sum_{n\leq x} d_{4}(n)  = C_1 x\log^3 x + C_2 x\log^2x + C_3 x\log x + C_4x + \vartheta(4.48 x^{3/4} \log x), \end{equation}
where
$$C_{1} = 1/6, \quad C_{2} = 0.654\ldots, \quad C_{3} = 0.981\ldots, \quad C_{4} = 0.272\ldots,$$
are exact constants given in (\ref{d4}). Furthermore, when $x\geq 193$ we have
\begin{equation}\label{sun}
\sum_{n\leq x} d_{4}(n) \leq \frac{1}{3} x \log ^{3} x.
\end{equation}
\end{thm}
\noindent It may be noted that the result in (\ref{chaise}) is a sharper bound than (\ref{hall}) and (\ref{lounge}) for all $x\geq 2$. Other bounds of the form (\ref{sun}) are possible: we have selected one that is nice and neat, and valid when $x$ is not too large. 

Sums of $d_{k}(n)$ can be used to obtain upper bounds on class numbers of number fields. Let $K$ be a number field of degree $n_{K}= [K:\mathbb{Q}]$ and discriminant $d_{K}$. Also, let $r_{1}$ (resp.\ $r_{2}$) denote the number of real (resp.\ complex) embeddings in $K$, so that $n_{K} = r_{1} + 2r_{2}$. Finally, let 
$$b= b_{K} = \left( \frac{n_{K}!}{n_{K}^{n_{K}}}\right) \left(\frac{4}{\pi}\right)^{r_{2}}|d_{K}|^{1/2},$$
denote the Minkowski bound, and let $h(K)$ denote the class number.
Lenstra \cite[\S 6]{Lenstra} --- see also Bordell\`{e}s \cite[Lem.\ 1]{Border} --- proved that
\begin{equation}\label{duck}
h_{K} \leq \sum_{m\leq b} d_{n_{K}}(m).
\end{equation}
We note that we need only upper bounds on (\ref{s2}) to give bounds on the class number.  Bordell\`{e}s used (\ref{hall}) in its weaker form
\begin{equation}\label{games}
\sum_{n\leq x} d_{k}(n) \leq 2x(\log x)^{k-1}, \quad (x\geq 6),
\end{equation}
to this end. The bound obtained on $h_{K}$ is not the sharpest possible for all degrees. For example, much more work has been done on quadratic extensions (see \cite{Le}, \cite{Louboutin}, and \cite{RamL}).  Using Theorem \ref{parlour} we are able to give an improved bound on the class number of quartic number fields.
\begin{cor}\label{chalk}
Let $K$ be a quartic number field with class number $h_{K}$ and Minkowski bound~$b$. Then, if $b\geq 193$ we have 
$$ h_{K} \leq \frac{1}{3} x\log^3x.$$
\end{cor}
We note that it should be possible to use Corollary \ref{chalk} to improve slightly Lemma 13 in \cite{Deb} when $n=4$.


Turning to the sum in (\ref{s2a}), we first state the following result by Ramanujan \cite{Ram2} 
\begin{equation}\label{bravo}
\sum_{n\leq x} d(n)^{2} - x(A \log^{3} x + B \log^{2} x + C\log x + D) \ll x^{3/5 + \epsilon}.
\end{equation}
The error in (\ref{bravo}) was improved by Wilson \cite{Wilson} to $x^{1/2 + \epsilon}$. The constants $A, B, C, D$ can (these days) be obtained via Perron's formula. Of note is an elementary result by Gowers \cite{Gowers}, namely, that
\begin{equation}\label{kitchen}
\sum_{n\leq x} d(n)^{2} \leq x(\log x + 1)^{3} \leq 2x (\log x)^{3}, \quad (x\geq 1).
\end{equation}
This is used by Kadiri, Lumely, and Ng \cite{Kadiri} in their work on zero-density estimates for the zeta-function. Although one would expect some lower-order terms, the bound in (\ref{kitchen}) is a factor of $2\pi^{2} \approx 19.7$ times the asymptotic bound in (\ref{s2a}), whence one should be optimistic about obtaining a saving. We obtain such a saving in our second main result.



\begin{thm}\label{billiard}
For $x\geq 2$ we have 
$$\sum_{n\leq x} d(n)^{2} = D_1 x\log^3 x + D_2 x\log^2x + D_3 x\log x + D_4x + \vartheta \left(9.73 x^\frac{3}{4} \log x + 0.73 x^\frac{1}{2} \right),$$ where
$$D_{1} = \frac{1}{\pi^2}, \quad D_{2} = 0.745 \ldots, \quad D_{3} = 0.824 \ldots, \quad D_{4} = 0.461\ldots,$$
are exact constants given in (\ref{dn2}). Furthermore, for $x\geq x_j$ we have $$\sum_{n\leq x} d(n)^2 \leq K x \log ^{3} x,$$ where one may take $\{K, x_j\}$ to be, among others, $\{\tfrac{1}{4}, 433\}$ or $\{1, 7\}$.
\end{thm}

The outline of this article is as follows. Theorem \ref{parlour} is proved in Section \ref{boat}. A similar process would give good explicit bounds on $\sum_{n\leq x} d_{k}(n)$. We have not pursued this, but the potential for doing so is discussed. We then present Theorem \ref{billiard} in Section \ref{ship}.

\section{Bounding $d_{4}(n)$}\label{boat}
Since $d_{4}(n) = d(n) * d(n)$, where $*$ denotes Dirichlet convolution, the hyperbola method gives us that
\begin{equation*}
\sum_{n\leq x} d_{4}(n) = 2 \sum_{a\leq \sqrt{x}} d(a) \sum_{n\leq x/a} d(n) - \left( \sum_{n\leq \sqrt{x}} d(n)\right)^{2}.
\end{equation*}

Using (\ref{study}), we arrive at 
\begin{equation}\label{daisy}
\begin{split}
\sum_{n\leq x} d_{4}(n) &= 2x \left[ \left(\log x + 2\gamma -1\right) S_{1} (\sqrt{x}) - S_{2}(\sqrt{x})\right] + 2\sum_{a\leq \sqrt{x}} d(a) \Delta \left(\frac{x}{a} \right) \\
&\quad - \left\{\sqrt{x}\left(\frac{1}{2} \log x + 2\gamma -1\right) + \Delta(\sqrt{x})\right\}^{2},
\end{split}
\end{equation}
where \begin{equation*}
S_{1}(x) = \sum_{n\leq x} \frac{d(n)}{n}, \text{ and } \quad S_{2}(x) = \sum_{n\leq x} \frac{d(n)\log n}{n}.
\end{equation*}
The absolute value of the sum on the right-hand side of (\ref{daisy}) can also be bounded above by
$$2 \alpha x^{1/2} S_{3}(\sqrt{x}),$$ where $$S_{3}(x) = \sum_{n\leq x} \frac{d(n)}{\sqrt{n}}.$$
We can approximate $S_1$, $S_2$, and $S_3$ with partial summation and the bound in (\ref{study}). We note that for applications in \S \ref{ship} we need only concern ourselves with values of $x\geq 1$.

Berkan\'{e}, Bordell\`{e}s, and Ramar\'{e} \cite[Cor.\ 2.2]{Berk},  give a bound for $S_{1}(x)$. As noted by Platt and Trudgian in \cite[\S 2.1]{Platt} their constant $1.16$ should be replaced by $1.641$ as in Riesel and Vaughan \cite[Lem.\ 1]{RV}. To obtain an error term in Theorem \ref{parlour} of size $x^{3/4} \log x$ we should like an error term in $S_{1}(x)$ of size $x^{-1/2}$, which is right at the limit of what is achievable.  We follow the method used by Riesel and Vaughan (\cite[pp. 48--50]{RV}) to write
\begin{equation}\label{chef}
S_{1}(x) = \frac{1}{2} \log^{2} x + 2 \gamma \log x + \gamma^{2} - 2 \gamma_{1} + \vartheta(c x^{-1/2}),
\end{equation}
with $c=1.001$ for $x\geq 6 \cdot 10^5$. One can directly check that this also holds for $2\leq x < 6 \cdot 10^5$. Taking larger values of $x$ reduces the constant $c$, but not to anything less than unity.

For $S_{2}$, we have
\begin{equation}\label{S2}
S_{2}(x) = \frac{1}{3} \log^{3} x + \gamma \log^{2} x + 2\gamma\gamma_{1} - \gamma_{2} + E_{2}(x),
\end{equation}
where
\begin{equation*}
E_{2}(x) = \frac{\Delta(x)\log x}{x} - \int_{x}^{\infty} \frac{(\log t -1)\Delta(t)}{t^{2}}\, dt.
\end{equation*}
Using the bound in (\ref{study}) with $\alpha=0.397$ and $x_0=5560$ we have
\begin{equation*}
|E_{2}(x)| \leq \alpha\left( 3 + \frac{2}{\log x_{0}}\right) x^{-1/2} \log x, \qquad (x\geq x_{0}).
\end{equation*}

Lastly, for $S_{3}$ we have
\begin{equation}\label{S3}
S_{3}(x) = 2x^{1/2} \log x + 4(\gamma -1)x^{1/2} + E_{3}(x),
\end{equation}
where
\begin{equation}\label{box}
E_{3}(x) = \frac{\Delta(x)}{x^{1/2}} - \frac{1}{2} \int_{x}^{\infty} \frac{\Delta(t)}{t^{3/2}}\, dt.
\end{equation}
For the integral in (\ref{box}) to converge we need to use a bound of the form $\Delta(t)\ll t^{1/2 - \delta}$. The only such explicit bound we know of is Theorem 1.2 in \cite{Berk}: as pointed out in \S \ref{intro} this improves on results only for large values of $x$. Instead, since $S_{3}(x)$ has a relatively small contribution to the total error, we can afford a slightly larger bound on $E_{3}(x)$. In writing the error as
\begin{equation*}
E_{3}(x) = 3 - 2\gamma + \frac{\Delta(x)}{x^{1/2}} + \frac{1}{2} \int_{1}^{x} \frac{\Delta(t)}{t^{3/2}}\, dt
\end{equation*}
we can apply the bound in (\ref{study}) and the triangle inequality to get
\begin{align*}
|E_{3}(x)| & \leq (3-2\gamma) + \alpha + \frac{1}{2} \alpha \log x \\
& \leq \log x\left(\frac{\alpha}{2} + \frac{3 - 2\gamma + \alpha}{\log x_{0}}\right):= \beta \log x, \qquad (x\geq x_{0}).
\end{align*}

Thus, the bounds in (\ref{chef}), (\ref{S2}), and (\ref{S3}) can be used in (\ref{daisy}) to prove Theorem \ref{parlour}
\begin{align} \label{d4}
\sum_{n\leq x} d_{4}(n) &= C_1 x\log^3 x + C_2 x\log^2x + C_3 x\log x + C_4x + E(x)
\end{align}
where $$|E(x)| \leq F_1 x^\frac{3}{4}\log x,$$ and we have
\begin{align*}
&C_1=\frac{1}{6}, \quad C_2= 2\gamma-\frac{1}{2} , \quad C_3= 6\gamma^2-4\gamma-4\gamma_1+1 , \\
& C_4= 4\gamma^3-6\gamma^2+ 4\gamma -12\gamma \gamma_1 + 4\gamma_1+2\gamma_2 -1 , \text{ and } F_1= 2c+6\alpha+ \frac{2\alpha}{\log x_{0}}.
\end{align*}

Recalling that $c=1.001, \alpha = 0.397$, and $x_{0} = 5560$, we prove the theorem for $x\geq 5560^2$. We directly calculated the partial sums of $d_4(n)$ to confirm that the bound in (\ref{chaise}) also holds for $2 \leq x < 5560^2$.

We could bound the partial sums of $d_{k}(n)$ by generalising the previous method. When $k$ is even, one can use $d_{k}(n) = d_{k-2}(n) * d(n)$, and when $k$ is odd, one can use $d_{k}(n) = d_{k-1}(n) * 1$. We have not pursued this, but for small values of $k$, this is likely to lead to decent bounds. We expect that the error term could, potential,  blow up on repeated applications of this process. One may also consider a more direct approach using a `hyperboloid' method, that is, considering $d_{3}(n) = 1 * 1 * 1.$ We have not considered this here.

\section{Bound on $\sum_{n\leq x} d(n)^2$}\label{ship}
We define $$H(s) = \sum_{n=1}^\infty \frac{h(n)}{n^{s}} = \frac{1}{\zeta(2s)},$$ and let $$H^{*}(s) = \sum_{n\geq 1} |h(n)| n^{-s} = \prod_{p} \left( 1 + \frac{1}{p^{2s}}\right) = \frac{\zeta(2s)}{\zeta(4s)},$$ whence both $H(s)$ and $H^{*}(s)$ converge for $\Re(s)> \frac{1}{2}$. Referring to (\ref{s1}), we therefore have $d(n)^2 = d_4(n) \ast h(n)$. Hence, we can write $$\sum_{n\leq x} d(n)^2 = \sum_{a\leq x} h(a) \sum_{b\leq \frac{x}{a}} d_4(b).$$

This leads to the bound
\begin{align} \label{dn2}
\sum_{n\leq x} d(n)^2 = D_1 x\log^3 x + D_2 x\log^2x + D_3 x\log x + D_4x + \sum_{a\leq x} h(a)E\left(\frac{x}{a}\right),
\end{align}
where
\begin{align*}
D_1&=C_1 H(1),\quad D_2=C_2H(1)+3C_1H'(1),\quad D_3=C_3H(1)+2C_2H'(1)+3C_1H''(1),  \\
D_4&=C_4H(1)+C_3H'(1)+C_2H''(1)+C_1H^{(3)}(1).
\end{align*}
Furthermore, we require a bound on $E(x)$ which holds for $x\geq 1$. Adapting the error term in (\ref{d4}) to hold for the desired range, we can write
$$\left| \sum_{a\leq x} h(a)E\left(\frac{x}{a}\right) \right| \leq F_1 H^{*}(\tfrac{3}{4}) x^\frac{3}{4} \log x + (1-C_4)x^\frac{1}{2}.$$
We also note the following exact values, for ease of further calculations:
\begin{equation}\label{tea}
\begin{split}
H(1) &= \frac{6}{\pi^{2}}, \quad
H'(1) = -\frac{72 \zeta'(2)}{\pi^{4}}, \quad
H''(1) = \frac{1728 \zeta'(2)^{2}}{\pi^{6}} - \frac{144 \zeta''(2)}{\pi^{4}},\\
H^{(3)}(1)&= -\frac{62208 \zeta'(2)^{3}}{\pi^{8}} + \frac{10368 \zeta'(2) \zeta''(2)}{\pi^{6}} - \frac{288 \zeta^{(3)}(2)}{\pi^{4}}.\\
\end{split}
\end{equation}

As an aside, our result in Theorem \ref{billiard} could have been achieved with Ramar\'{e}'s general result in \cite[Lem.\ 3.2]{Ram}, as modified in \cite[Lem.\ 14]{TQ} (with the constants repaired as in \cite{RV}), but some further generalisation would have been necessary. Instead, we proceeded directly as above.

\subsection*{Acknowledgements}
We wish to thank Olivier Ramar\'{e} for many rollicking conversations on this topic.

 \end{document}